\documentclass[11pt]{amsart}
\usepackage{latexsym,amsmath,amsthm,amssymb}
\usepackage{pslatex}
\usepackage{color}
\usepackage[colorlinks=true,linkcolor=blue,citecolor=nicegr]{hyperref}
\definecolor{nicegr}{rgb}{0,.5,0}

\title[Primitive Central Idempotents]{Primitive Central Idempotents \\ of the Group Algebra}
\author{Robin Endelman}
\address{University College of the Fraser Valley,
Abbotsford, BC  V2S 7M8, Canada}
\email{robin.endelman@ucfv.ca}
\author{Manash Mukherjee}
\address{University College of the Fraser Valley,
Abbotsford, BC  V2S 7M8, Canada}
\email{manash.mukherjee@ucfv.ca}

\keywords{primitive idempotent, group algebra, irreducible representation}

\subjclass[2000]{20C05, 16S34}

\newtheorem{thm}{Theorem}[section]
\newtheorem{lem}[thm]{Lemma}
\newtheorem{prop}[thm]{Proposition}
\newtheorem{cor}[thm]{Corollary}
\newtheorem{fact}[thm]{Fact}

\theoremstyle{definition}
\newtheorem{defn}[thm]{Definition}

\theoremstyle{remark}

\newtheorem{rem}{Remark}

\numberwithin{equation}{section}

\newcommand{\Z}{\mathbb{Z}}
\newcommand{\C}{\mathbb{C}}

\newcommand{\CG}{\C[G]}
\newcommand{\CN}{\C[N]}
\newcommand{\Zcg}{Z(\CG)}

\newcommand{\Oa}{{\mathcal O}_\alpha}

\begin{document}

\begin{abstract}
An approach to representations of finite groups is presented without
recourse to character theory. Considering the group algebra $\CG$ as an algebra of linear maps on $\CG$ (by left multiplication), we derive the primitive central idempotents as a simultaneous eigenbasis of the centre, $\Zcg$.   We apply this framework to obtain the irreducible representations of a class of finite meta-abelian groups. In particular, we give a general construction of the isomorphism between simple blocks of $\CG$ and the corresponding matrix algebra where $G$ can be any finite group.

\end{abstract}

\maketitle

\tableofcontents

\section{Introduction}

In this paper, we present an approach to determining the decomposition of the group algebra, $\CG$, of a finite group $G$, over the field of complex numbers, $\C$.  We work solely within the group algebra without recourse to the character theory of $G$, obtaining the irreducible representations of $G$.

There are, of course, many ways of approaching the representations of finite groups \cite{Serre, GW, CR-old, L}, however, our starting point is the centre of the group algebra.  The very concept of primitive central idempotents (along with their properties) arises naturally as a consequence of the diagonalizability of $\Zcg$.

We consider the group algebra acting on itself by left multiplication (that is, considering the regular representation of the algebra $\CG$),
$$\CG \times \CG \rightarrow \CG : (x,y) \mapsto xy$$
and derive the primitive central idempotents of the group algebra as a simultaneous eigenbasis of its centre, $\Zcg$ (Theorem \ref{thm:ebasis}, Section \ref{sec:J}).  An immediate consequence, among others, of this framework is that the two-sided ideal generated by a primitive central idempotent has a one-dimensional centre.

Our point of view is illustrated first for abelian groups in Section \ref{sec:abel}. In Section \ref{sec:meta}, we apply this framework to a large class of nonabelian groups, namely (meta-abelian) semi-direct products, $G = N \rtimes H$ where both $N$ and $H$ are abelian. For this class of groups, we provide a general formula for the primitive central idempotents of $\CG$ (Section \ref{sec:pci}). To obtain the isomorphism between each simple block of $\CG$ and the corresponding matrix algebra, we give a general construction that is applicable for all finite non-abelian groups (Proposition \ref{prop:iso}).  Furthermore, a formula is given for counting the number of conjugacy classes in $N \rtimes H$ (Corollary \ref{cor:Jcount}), which also holds when $H$ is nonabelian (Remark \ref{rem:nonabel}).

Throughout this paper $G$ denotes a finite group.  

\subsection{Acknowledgements}

We  would like to thank  Stephanie van Willigenburg (Department of Mathematics, University of British Columbia, Canada) for her comments and careful reading of the manuscript.  R.E.   is grateful to the Faculty of Research and Graduate Studies (UCFV) for continued support. 

\section{Primitive central idempotents of \texorpdfstring{$\CG$}{\CG}} 
\label{sec:J}

In this section, we obtain the primitive central idempotents of the group algebra as a simultaneous eigenbasis for the centre.  

\begin{defn}\label{defn:idempt}

(i) An element $u$ in an algebra $A$ is {\it idempotent} if $u^2 = u$.  
Two idempotents $u_1$ and $u_2$ are {\it orthogonal} if 
$$u_1 u_2 = u_2 u_1 = 0.$$  

(ii) An idempotent $u$ is called {\it primitive} if $u$ cannot be written as a sum, $u = u_1 + u_2$, where $u_1$ and $u_2$ are orthogonal idempotents.  

(iii) $J$ is called a {\it primitive central idempotent} if $J$ belongs to the centre of $A$, $J$ is idempotent, and $J$ cannot be written as a sum, $$J = I_1 + I_2,$$ 
where $I_1$ and $I_2$ are orthogonal central idempotents.

\end{defn}

We first show that the centre, $\Zcg$, consists of diagonalizable elements, and then obtain the primitive central idempotents of $\CG$ as the simultaneous eigenbasis of the centre.  

\noindent {\sf Notation}: For $x \in \CG$, $\bar{x}$ will denote the complex conjugate of $x$, and $x^\dagger$ the conjugate transpose of $x$, where $x$ is considered as a linear operator on $\CG$ by left multiplication.

\begin{prop} \label{prop:diag}
Every element of $\Zcg$ is diagonalizable (operating by left multiplication on $\Zcg$).
\end{prop}

\begin{proof} 
As an operator on $\CG$ by left multiplication, each group element $g$ is a permutation matrix, and therefore $g^{\dagger} = g^{-1}$.  Thus, for any $x = \sum_{g \in G} c_g g \in \CG$, $x^\dagger = \sum_{g \in G} \bar{c}_g g^{-1} \in \CG$.  If $z \in \Zcg$, $z z^\dagger = z^\dagger z$, and consequently, $z$ is diagonalizable on $\CG$.  Now the proposition follows from the fact that the centre, $\Zcg$, is an invariant subspace for $z$.
\end{proof}

\begin{rem}
If $z \in \Zcg$, then $z^\dagger \in \Zcg$.  
In particular, the standard basis elements, $\Omega_g = \sum_{x \in K_g} x$, of $\Zcg$, satisfy $\Omega_g^{\dagger} = \Omega_{g^{-1}}$.  [Here, $K_g$ denotes the conjugacy class of $g$.] 
\end{rem}

By Proposition \ref{prop:diag}, any basis of $\Zcg$ is a commuting set of diagonalizable linear operators on $\Zcg$ (acting by left multiplication), and hence, there is a simultaneous eigenbasis of 
$\Zcg$.

\begin{thm} \label{thm:ebasis} 
There is a unique simultaneous eigenbasis, $\{J_p : 1 \leq p \leq k\}$, of $\Zcg$ such that:
\begin{itemize}
\item[${\rm (i)}$] 
$J_p^2 = J_p$
\item[${\rm (ii)}$] $J_p J_q = 0$ for $p \neq q$, 
\item[${\rm (iii)}$] the group identity, $e = \sum_p J_p$, and 
\item[${\rm (iv)}$] $J_p$ is a primitive central idempotent.
\end{itemize}
\end{thm} 

\begin{proof} 
(i) and (ii):  Let $\{\tilde{J_p} : 1 \leq p \leq k\}$ be a simultaneous eigenbasis of $\Zcg$.
For any $p, q$, viewing $\tilde{J_q}$ as an eigenvector of $\tilde{J_p}$, and $\tilde{J_p}$ as an eigenvector of $\tilde{J_q}$, we have 
\begin{eqnarray*}
\tilde{J_p} \tilde{J_q} & = & \lambda \tilde{J_q} \\
\tilde{J_q} \tilde{J_p} & = & \mu \tilde{J_p}
\end{eqnarray*}
for some $\lambda, \mu \in \C$.  Since $\tilde{J_p}$ and $\tilde{J_q}$ commute, $\lambda \tilde{J_q} = \mu \tilde{J_p}$.  
Hence, by linear independence of $\tilde{J_p}$ and $\tilde{J_q}$, we have $\lambda = \mu = 0$ for $p \neq q$, and for $p = q$, $\tilde{J_p}^2 = \lambda \tilde{J_p}$.  
Now, each $\tilde{J_p}$ is scaled to obtain the unique eigenvector $J_p$ satisfying $J_p^2 = J_p$, and $J_p J_q = 0$.  

\noindent (iii)  Set $e = \sum_q c_q J_q$.  Then, for each $p$, 
$J_p = e J_p = c_p J_p$, and hence, $c_p = 1$ for all $p$.

\noindent (iv)  If $J_p$ decomposes as a sum, $J_p = I_1 + I_2$, of two mutually orthogonal central idempotents $I_1$ and $I_2$, 
then, $\{I_1, I_2\} \cup \{J_q | q \neq p\}$ would be $k+1$ linearly independent central elements, exceeding the dimension of $\Zcg$.  
\end{proof}

\begin{rem}
Since $J_p$ is central, $A_p = (\CG) J_p$ is a two-sided ideal of $\CG$.  By Theorem \ref{thm:ebasis}, 
it follows that the group algebra $\CG$ is a direct sum of the two-sided ideals $\{A_p\}$,
$$\CG = \bigoplus_{p=1}^k A_p.$$
\end{rem}

\begin{rem}\label{rem:Z1d}
For each $p$, $J_p$ spans the 1-dimensional centre of $A_p$:  if $a \in Z(A_p)$, then $a \in \Zcg$, and hence, $a = a J_p = \lambda J_p$, where the first equality follows from $a \in A_p$ and the second equality follows since $J_p$ is an eigenvector of $a \in \Zcg$.   
\end{rem}

\begin{rem}
It follows from Remark \ref{rem:Z1d} that $A_p$ cannot be decomposed as a direct sum of two-sided ideals. 
\end{rem}
 
It is well known that the blocks $A_p$ are simple, and therefore isomorphic to matrix algebras (\cite{L}, for example).

\section{Abelian groups} 
\label{sec:abel}

We give several illustrations of our point of view in Theorem \ref{thm:ebasis}, that 
the primitive central idempotents are a simultaneous eigenbasis of $\Zcg$.  

\subsection{Irreducible representations}

We begin with the well-known result:  A finite group $G$ is abelian if and only if $\CG \cong \C^k$.

\begin{proof}  If $G$ is abelian, then $\CG$ is a commutative algebra.  Let $\{J_p\}$ be the simultaneous eigenbasis of $\CG$, as in Theorem \ref{thm:ebasis}.  
Then, for every $g \in G$, $g J_p = \lambda_p(g) J_p$ where $\lambda_p(g) \in \C$.  
Hence, $\CG = \oplus_{p=1}^k \C J_p \cong \C^k$.  The converse is trivial.  
\end{proof}

\begin{rem}
The ideals $\C J_p$ provide the one-dimensional (irreducible) representations of $G$.  Since $G$ is finite, $|\lambda_p(g)| = 1$ for all $g \in G$.  (The functions $\lambda_p$ are the characters of $G$.)
\end{rem}

\subsection{Cyclic groups} 

We illustrate the simultaneous eigenbasis for a cyclic group $G = \Z_k = \{e, r, r^2, ..., r^{k-1}\}$, with generator $r$ of order $k$.  

\begin{prop} \label{prop:cyclic}
Let $G = \langle r \rangle$ be a cyclic group of order $k$ with generator $r$, and $\omega = \exp(i2\pi/k)$. For $1 \leq p \leq k$, orthogonal idempotents 
$$J_p = \frac{1}{k} \sum_{m=1}^{k} \omega^{-mp} r^m$$
form the simultaneous eigenbasis of $\CG$. 
\end{prop} 

\begin{proof}
As the minimal polynomial of $r$ is $x^k - 1$, the distinct eigenvalues of $r$ are $\lambda = \omega^p$, $1 \leq p \le k$, where $\omega = \exp(i2\pi/k)$.  Let $v = \sum_{m=1}^{k} a_m r^m$ be an eigenvector of $r$ corresponding to eigenvalue $\lambda$.
Then $r v = \lambda v$ implies 
$$a_k = \lambda a_1 = \lambda^2 a_2 = \cdots = \lambda^{k-1} a_{k-1}.$$ 
Taking $a_k = 1$, the eigenvectors of $r$ with eigenvalue $\lambda = \omega^p$ are 
multiples of 
$v_p = \sum_{m=1}^{k} \omega^{-mp} r^m$.
Now, 
\begin{eqnarray*}
v_p v_q &=& \sum_{m=1}^{k} \omega^{-mp} r^m v_q \\
 &=& \sum_{m=1}^{k} \omega^{-m(p-q)} v_q 
 = \left\{ \begin{array}{ll}
               0 & q \neq p \\
               k v_q & q = p
             \end{array} \right.
\end{eqnarray*}
Thus the desired (orthogonal) idempotents are $J_p = \frac{1}{k} v_p$.
\end{proof}

\subsection{Arbitrary finite abelian groups} 

For the group algebra of an arbitrary finite abelian group, the simultaneous eigenbasis of Theorem \ref{thm:ebasis} is obtained by Proposition \ref{prop:cyclic} above, together with the following result:

\begin{prop}\label{prop:product}
For any finite groups, $G_i$, 
$$\C[G_1 \times G_2 \times \cdots \times G_n] \simeq \C[G_1] \otimes \C[G_2] \otimes \cdots \otimes \C[G_n]$$
\end{prop}

\begin{proof}
The isomorphism is given by defining 
$$\phi(g_1, g_2, ..., g_n) = g_1 \otimes g_2 \otimes \cdots \otimes g_n,$$ and extending linearly. 
\end{proof}

For an arbitrary finite abelian group $G = \Z_{n_1} \times \cdots \times \Z_{n_k}$, where $\Z_{n_i}$ is the cyclic group of order $n_i$, the idempotent eigenbasis for each $\C[\Z_{n_i}]$ is given by Proposition \ref{prop:cyclic}.  Let $J^{(i)}_p$, $1 \leq p \leq n_i$ denote these idempotents for each $\Z_{n_i}$. Then the idempotent eigenbasis of $\CG$ is given by 
$$J^{(1)}_{p_1} J^{(2)}_{p_2} \cdots J^{(k)}_{p_k}$$
where $1 \leq p_i \leq n_i$ for each $i$.

\begin{rem}
For $G$ abelian, the primitive central idempotents of $\CG$ are primitive idempotents.  When $G$ is nonabelian, $\Zcg$ is a proper subalgebra of $\CG$.  In this case, each primitive central idempotent decomposes as a sum of orthogonal primitive idempotents, 
$$J_p = u_1^p + u_2^p + \cdots + u_{n_p}^p$$ 
[See the Appendix (Section \ref{sec:primu}) for characterizations of primitive idempotents.]  
\end{rem}

In the following section, we apply the results of Sections \ref{sec:J} and \ref{sec:abel} to  a large and important class of nonabelian groups --- the semi-direct product of two abelian groups.  We construct the irreducible representations without recourse to character theory or induced representations. 

\section{Meta-abelian semi-direct products} 
\label{sec:meta}

We consider a semi-direct product $G = N \rtimes H$, where $N$ and $H$ are finite abelian groups. We construct the primitive idempotents of $\CG$ in terms of the primitive (central) idempotents of $N$ and those of $H$.  Consequently, we obtain each primitive central idempotent of $\CG$ as a sum of primitive idempotents.  This leads to the explicit isomorphism of $\CG$ with the direct sum of matrix algebras.  

Let $P_N = \{v_a : 1 \leq a \leq |N| \}$ be the orthogonal idempotent basis (of Theorem \ref{thm:ebasis}) for $\C[N]$.  As $N$ is abelian, these idempotents are obtained by Propositions \ref{prop:cyclic} and \ref{prop:product}.  Furthermore, conjugation by $H$ is an algebra isomorphism on $\CN$, and hence, by the uniqueness of the idempotent eigenbasis, $P_N$ is an $H$-set.

\subsection{Primitive idempotents of \texorpdfstring{$\CG$}{\CG}} 
\label{sec:prim}

Let $\Oa = \{v_1, v_2, ..., v_{c_\alpha}\}$ be an orbit in $P_N$.  We set $\{u_p : 1 \leq p \leq |H_\alpha|\}$ to be the idempotent eigenbasis of $\C[H_\alpha]$, where $H_\alpha$ denotes the common stabilizer of any $v_i \in \Oa$.

\begin{prop} \label{prop:vuprim}
For each $v_i \in \Oa$, $1 \leq i \leq |\Oa|$, and each  $u_p \in \C[H_\alpha]$, $1 \leq p \leq |H_\alpha|$, 
the elements $v_i u_p$ are primitive idempotents in $\CG$.
\end{prop}

\begin{proof}
As $v_i$ and $u_p$ commute, $(v_i u_p)^2 = v_i u_p$.  
For primitivity, we use the fact (see Appendix, Fact \ref{fact:uprim}) that an idempotent $u \in \CG$ is primitive if and only if $u (\CG) u = \C u$.  
For $nh \in G$, 
$$v_i u_p (nh) v_i u_p 
= \lambda (v_i u_p~h~v_i u_p), \mbox{ for some $\lambda \in \C$,}$$
since $v_i$ is a simultaneous eigenvector of $\C[N]$.
If $h \notin H_\alpha$, $h v_i h^{-1} = v_j$ for some $j \neq i$ and the right hand side of the above expression vanishes.  
For $h \in H_\alpha$, $h u_p = \sigma u_p$ for some $\sigma \in \C$ and hence,
$$v_i u_p(nh)v_i u_p 
= \lambda \sigma (v_i u_p).$$
\end{proof}

\begin{prop}\label{prop:vuequiv}
For $v_1, v_2 \in \Oa$, and $1 \leq p, q \leq |H_\alpha|$, the primitive idempotents
$v_1 u_p$ and $v_2 u_q$ are equivalent if and only if $p = q$.
\end{prop}

\begin{proof}  

For $v_1, v_2 \in \Oa$, there is an $h \in H$ such that $v_1 = h v_2 h^{-1}$.  It follows that 
$h(v_1 u_p) h^{-1} = v_2 u_p$, and hence $v_1 u_p$ and $v_2 u_p$ are equivalent (see Appendix, Fact \ref{fact:xux}).

If $v_1 u_p$ and $v_2 u_q$ are equivalent, then (see Appendix, Fact \ref{fact:uxv}) there is an element $nh \in G$ so that $v_1 u_p (nh) v_2 u_q \neq 0$.  The left hand side is equal to $\lambda v_1 u_p u_q h v_2$, and is nonzero only when $p = q$. 
\end{proof}

\subsection{Primitive central idempotents of \texorpdfstring{$\CG$}{\CG}} 
\label{sec:pci}

For each $\Oa$, we sum the equivalent primitive idempotents of Proposition \ref{prop:vuequiv} in order to construct the primitive central idempotent basis for $\Zcg$.  Given an orbit $\Oa$, we define
$$
J^{[\alpha]}_p = \sum_{i=1}^{|\Oa|} v_i u_p
$$
for each $u_p \in \C[H_\alpha]$, $1 \leq p \leq |H_\alpha|$.

\begin{lem}\label{lem:Jprop}
The elements $\{J^{[\alpha]}_p\}$ 
are orthogonal central idempotents of $\CG$.  
Furthermore, $\displaystyle{\sum_{\alpha, p} J^{[\alpha]}_p = e}$.
\end{lem}

\begin{proof}
(a) $J^{[\alpha]}_p$ are idempotents:  For $v_1, v_2 \in \Oa$, the idempotents $v_1 u_p$ and $v_2 u_p$ are orthogonal.  Thus $J^{[\alpha]}_p$ is a sum of orthogonal idempotents, and hence is  idempotent.

(b) $J^{[\alpha]}_p$ are central:  We write $J^{[\alpha]}_p = S_\alpha u_p$, where $S_\alpha = \sum_{i=1}^{|\Oa|} v_i$.  
Then, for any $h \in H$, $h S_\alpha h^{-1} = S_\alpha$, and hence $h J^{[\alpha]}_p h^{-1} = J^{[\alpha]}_p$.  For $n \in N$ and $v_i \in \Oa$, $n v_i = \lambda_i v_i$ where $\lambda_i \in \C^\times$ (and $|\lambda_i| = 1$).  Thus, $n S_\alpha = \sum_{i} \lambda_i v_i$ commutes with $u_p$, so that $nJ^{[\alpha]}_p n^{-1} = J^{[\alpha]}_p$.

(c) $J^{[\alpha]}_p$ are orthogonal: For $\alpha \neq \beta$, $S_\alpha S_\beta = 0$ by orthogonality of the primitive idempotents of $\C[N]$, and hence, 
$J^{[\alpha]}_p J^{[\beta]}_q  = 0$.  
For $\alpha = \beta$ and $p \neq q$, $J^{[\alpha]}_p J^{[\alpha]}_q = 0$, by orthogonality of the primitive idempotents of $\C[H_\alpha]$.

(d) For each $\alpha$, $\sum_p u^{[\alpha]}_p = e$, where we denote by $u^{[\alpha]}_p$ the primitive idempotents of 
$\C[H_\alpha]$.  Thus, 
\begin{eqnarray*}
\sum_{\alpha,p} J^{[\alpha]}_p & = & \sum_{\alpha, p} S_\alpha u^{[\alpha]}_p = \sum_\alpha S_\alpha \\
 & = & \sum_{v \in P_N} v  = e.
\end{eqnarray*}
\end{proof}

\begin{thm}\label{thm:Jbasis}
The elements $\{J^{[\alpha]}_p\}$ are the 
primitive central idempotents of $\CG$.  
\end{thm}

\begin{proof}
The primitive summands, $v_i u_p$, of $J^{[\alpha]}_p$ are equivalent, and hence, by Fact \ref{fact:equivu}, they lie in the same block $A = (\CG)J$ of $\CG$ (where $J$ is a primitive central idempotent of $\CG$).  Thus, $J^{[\alpha]}_p \in A \cap Z(\CG) = \C J$.  
It follows that $J^{[\alpha]}_p = J$.  By Lemma \ref{lem:Jprop}(d) above, these are all the primitive central idempotents.
\end{proof}

\begin{cor}   \label{cor:Jcount}
Let $G = N \rtimes H$, where $N$ and $H$ are abelian.  With the stability subgroups, $H_\alpha$, defined in Section \ref{sec:prim}, 
$$\dim(\Zcg) = \sum_\alpha |H_\alpha|$$
where the sum is taken over $\alpha$ corresponding to distinct orbits $\Oa$.  
Hence, the number of conjugacy classes of $G$ is equal to the sum 
of the number of conjugacy classes of the stability subgroups.  
\end{cor}

\begin{proof}
By Theorem \ref{thm:Jbasis}, $\{J^{[\alpha]}_p\}$ forms a basis for $\Zcg$.  For each $\alpha$ corresponding to an orbit $\Oa$, the number of $J^{[\alpha]}_p$ (which is equal to the number of $u^{[\alpha]}_p$), is $|H_\alpha|$.  Hence, 
$\dim(\Zcg) = \sum_\alpha |H_\alpha|,$
where ($H_\alpha$ being abelian) $|H_\alpha|$ is the number of conjugacy classes of $H_\alpha$. 
\end{proof}

\subsection{Irreducible representations of \texorpdfstring{$\CG$}{\CG}} 
\label{sec:irrep}

Let $A = (\CG)J$ be a simple block generated by a primitive central idempotent, $J=\sum_{i=1}^{m} u_i$, where the primitive idempotents 
$u_i$ determine the decomposition of $A$ into isomorphic minimal left ideals (as in the Appendix, Section \ref{sec:primu}),  
$A = \bigoplus_{j=1}^{m} A u_j$. Each ideal $A u_j$ is of dimension $m$, and is a sum of one-dimensional vector spaces, $A u_j = \oplus_{i=1}^{m} u_i A u_j$.
A basis element $e_{ij} \in u_i A u_j$ is determined by an element $a \in A$ for which $u_i a u_j \neq 0$.  The idempotents $u_i$ and $u_j$ being equivalent, we can find an invertible element $x \in A$ satisfying $x u_j x^{-1} = u_i$.  Then $u_i x u_j = x u_j \neq 0$.

\begin{prop} \label{prop:iso}
Let $A = (\CG)J$ where $J= \sum_{i=1}^m u_i$ is a primitive central idempotent of $\CG$ and $u_i$ are orthogonal primitive idempotents.  For each $1 \leq i \leq m$, let $x_i \in A$ be such that $x_i u_1 x_i^{-1} = u_i$ and 
$$e_{ij} = u_i (x_i x_j^{-1}) u_j.$$
Then, $e_{ij}$ defines an isomorphism of $A$ with ${\mathbb M}_m(\C)$. 
\end{prop}

\begin{proof}
\begin{eqnarray*}
e_{ij} e_{jk} & = & u_i x_i x_j^{-1} u_j x_j x_k^{-1} u_k \\
  & = & u_i x_i u_1 x_k^{-1} u_k = u_i x_i x_k^{-1} u_k \\
  & = & e_{ik} 
\end{eqnarray*}
For $j\neq k$, $u_j u_k =0$ and hence, $e_{ij} e_{kl} = 0$.
Mapping $e_{ij}$ to the standard basis elements $E_{ij}$ of the matrix algebra ${\mathbb M}_{m}(\C)$, defines an isomorphism of 
$A$ with ${\mathbb M}_m(\C)$.
\end{proof}

By this construction, we obtain the irreducible representation for each simple summand of $\CG$.

We now apply Proposition \ref{prop:iso} to the meta-abelian semi-direct product, $G = N \rtimes H$.
For fixed $\alpha$ and $p$, let us consider a simple block, $A = (\CG)J^{[\alpha]}_p$, generated by a primitive central idempotent 
$J^{[\alpha]}_p = \sum_{i=1}^{c_\alpha} v_i u_p$, and   
set $w_i \equiv v_i u_p$.    
Each ideal $A w_j$ is of dimension $c_\alpha = |\Oa| = [H:H_\alpha]$, and is a sum of one-dimensional vector spaces, $A w_j = \oplus_{i=1}^{c_\alpha} w_i A w_j$.  

To define a basis element, $e_{ij} \in w_i A w_j$,  it is sufficient to consider $a = nh \in G$, and furthermore, since $v_{i}$ (and hence $w_i$) are simultaneous eigenvectors of $N$, it suffices to consider elements $a \in H$.

Now, we choose an element $h_i$ from each coset of $H_\alpha$ in $H$,  set $v_{i} = h_i v_{1} h_i^{-1}$ (and hence $w_i = h_i w_1 h_i^{-1}$), and by Proposition \ref{prop:iso}, define 
$$e_{ij} = w_i h_i h_j^{-1} w_j \in w_i A w_j.$$
Then, $e_{ij} e_{kl} = 0$, $j\neq k$ and $e_{ij} e_{jk}=e_{ik}$.   With these $e_{ij}$, we have the required isomorphism of $A$ with the matrix algebra ${\mathbb M}_{c_\alpha}(\C)$.

\subsection{Example}
\label{sec:exm}

Let $G = A_4 = N \rtimes H$ with $N = \Z_2 \times \Z_2 = \langle x \rangle \times \langle y \rangle = \{e, x, y, z\}$, where $z = xy = yx$, $x^2 = y^2 = e$, and $H = \Z_3 = \langle t \rangle$, where $t^3 = e$.  The action of $H$ on $N$ is defined by 
$$txt^{-1} = z ~,~ tyt^{-1} = x ~,~ tzt^{-1} = y.$$

\subsubsection{The decomposition of $\C[A_4]$}  
\label{sec:dcmp}

The primitive central idempotents of $\C[N]$ are obtained by Propositions \ref{prop:cyclic} and \ref{prop:product}: 
\begin{eqnarray*}
v_0 & = & [(e+x)/2][(e+y)/2] = (e + x + y + z)/4 \\
v_1 & = & [(e+x)/2][(e-y)/2] = (e + x - y - z)/4 \\
v_2 & = & [(e-x)/2][(e+y)/2] = (e - x + y - z)/4 \\
v_3 & = & [(e-x)/2][(e-y)/2] = (e - x - y + z)/4
\end{eqnarray*}
The action of $H$ on $P_N = \{v_0, v_1, v_2, v_3\}$ determines two orbits, $${\mathcal O}_0  =  \{v_0\} ~~,~~ {\mathcal O}_1 = \{v_1, v_2, v_3\}$$  
with corresponding stability subgroups 
$$H_0 = H = \{e, t, t^2\} ~~,~~ H_1 = \{e\}.$$ 
The primitive central idempotents of $H_0$ are also given by Proposition \ref{prop:cyclic}:
$$
u_0  =  \frac{1}{3}(e + t + t^2) ~,~
u_1  =  \frac{1}{3}(e + \omega^2 t + \omega t^2) ~,~
u_2  =  \frac{1}{3}(e + \omega t + \omega^2 t^2).
$$

By Theorem \ref{thm:Jbasis}, the primitive central idempotents of $\C[A_4]$, as sums of primitive idempotents are:
\begin{eqnarray*}
J_a  & = &  v_0 u_a, \mbox{ for $a = 0, 1, 2$},\\
J_3  & = & v_1 + v_2 + v_3.
\end{eqnarray*}
We see that $(\CG)J_0$, $(\CG)J_1$, and $(\CG)J_2$ are of dimension 1, $(\CG)J_3$ is of dimension $|{\mathcal O}_1|^2 = 3^2$, and 
$$\C[A_4] \simeq \C \oplus \C \oplus \C \oplus {\mathbb M}_3(\C).$$ 

\subsubsection{A basis for $(\CG)J_3$}
\label{sec:bs}

In $J_3$, $w_i = v_i e = v_i$.  The action of $H$ on $P_N$ is given by
$$t^2 v_1 t^{-2} = v_2  ~,~   t v_1 t^{-1} = v_3 $$
and the coset representatives for $H_1$ are $h_1 = e$, $h_2 = t^2$, and $h_3 = t$.  
Thus, as in Section \ref{sec:irrep}, $e_{ij} = w_i h_i h_j^{-1} w_j$ gives the column spaces $Av_j$:
$$\left(
\begin{array}{r}
e_{11} \\
e_{21} \\
e_{31} 
\end{array} \right) = 
\left(
\begin{array}{r}
v_1 e v_1 \\
v_2 t^2 v_1 \\
v_3 t v_1 
\end{array} \right) = 
\left(
\begin{array}{r}
v_1 \\
t^2 v_1 \\
t v_1
\end{array} \right)$$
$$\left(
\begin{array}{r}
e_{12} \\
e_{22} \\
e_{32} 
\end{array} \right) = 
\left(
\begin{array}{c}
v_1 t v_2 \\
v_2  \\
v_3 t^2 v_2 
\end{array} \right) = 
\left(
\begin{array}{c}
t v_2 \\
v_2 \\
t^2 v_2
\end{array} \right)$$
and
$$\left(
\begin{array}{r}
e_{13} \\
e_{23} \\
e_{33} 
\end{array} \right) = 
\left(
\begin{array}{c}
v_1 t^2 v_3 \\
v_2 t v_3 \\
v_3 
\end{array} \right) = 
\left(
\begin{array}{c}
t^2 v_3 \\
t v_3 \\
v_3
\end{array} \right).$$
Thus, 
\begin{eqnarray*}
xJ_3 &=& v_1 - v_2 - v_3 = e_{11} - e_{22} - e_{33} \\ 
yJ_3 &=& -v_1 + v_2 - v_3 = -e_{11} + e_{22} - e_{33}\\
tJ_3 &=& t v_1 + t v_2 + t v_3 = e_{31} + e_{12} + e_{23}
\end{eqnarray*}
determine the three dimensional irreducible representation of $A_4$.

\begin{rem} \label{rem:nonabel}
If $G = N \rtimes H$ where $H$ is nonabelian (keeping $N$ abelian),
construction of the primitive central idempotents proceeds in a similar way as the meta-abelian case.  
In particular, Corollary \ref{cor:Jcount} also holds in this case:
$$\dim(\Zcg) = \sum_\alpha \dim(Z(\C[H_\alpha])),$$
where $H_\alpha$ represents the isomorphism class of the stability subgroups for the orbit $\Oa$. 
\end{rem}

\section{Appendix: Primitive idempotents}
\label{sec:primu}

An idempotent $u$ of an associative algebra $A$ is called {\it primitive} if $u$ cannot be written as a sum, $u = u_1 + u_2$, where $u_1^2 = u_1$, $u_2^2 = u_2$, and $u_1 u_2 = u_2 u_1 = 0$.  

We list several properties of primitive idempotents \cite{P} in a complex semisimple algebra, $A = \oplus_i A_i$, where $A_i $ are the simple summands.  

\begin{rem}
By standard arguments (among the simplest, see \cite{L}), the simple components $A_i$ are isomorphic to matrix algebras, ${\mathbb M}_{n_i}(\C)$.
\end{rem}

\begin{fact} \label{fact:uprim}
An idempotent $u$ is primitive if and only if $Au$ is a minimal left ideal; equivalently, $uAu = \C u$.  Hence, 
a primitive idempotent $u \in A$ belongs to a simple block $A_i$.  If $u \in A_i$, then $Au = A_i u$.   
\end{fact}

\begin{fact}\label{fact:uxv}
Minimal left ideals $Au_1$ and $Au_2$ are isomorphic if and only if $u_1 x u_2 \neq 0$ for some $x \in A$.  
\end{fact}

\begin{fact} \label{fact:xux}
Minimal left ideals $Au_1$ and $Au_2$ are isomorphic if and only if there exists an invertible element $x \in A$ such that $x u_1 x^{-1} = u_2$.
\end{fact}

From the above Facts it follows that
\begin{fact} \label{fact:equivu}
Given two primitive idempotents $u_1$ and $u_2$, $Au_1$ and $Au_2$ are isomorphic left ideals if and only if $u_1$ and $u_2$ belong to the same simple summand, $A_i$, of $A$.  
\end{fact}
\noindent 
In this case, we call $u_1$ and $u_2$ {\it equivalent} primitive idempotents.  

We fix $i$ and set $n = n_i$ and $B = A_i$.  
Let $e_{ij}$ be the elements of $B$ that correspond to $E_{ij}$ -- the standard basis elements of ${\mathbb M}_n(\C)$ with 1 in the $(i,j)$ position and 0 elsewhere.  If $J$ is the identity element of $B$, then the primitive idempotents $u_i = e_{ii}$ are characterized by

\begin{itemize}
\item $J = u_1 + \cdots + u_n$
\item $u_i^2 = u_i$
\item $u_i u_j = 0$ for $i \neq j$
\end{itemize}
Furthermore, the vector spaces $u_i B u_j$ are all one-dimensional: they are equal to $\C e_{ij}$ for $i \neq j$, and to $\C u_i$ when $i = j$.  
It follows that $B$ decomposes as a direct sum of isomorphic minimal left ideals
$$B = \oplus B u_i$$ 
where the dimension of each left ideal $B u_i$ is $n$.

\end{document}